\documentclass[11pt, a4paper]{article}

\usepackage{epsfig}
\usepackage{amssymb}
\usepackage{amsthm}
\usepackage{amsmath}
\usepackage{float}
\usepackage{multicol}
\usepackage{multirow}
\usepackage{cite}
\usepackage{enumerate}
\usepackage{xcolor}
\usepackage{tabu}

\def\ms{\medskip}
\def\nt{\noindent}

\definecolor{vividviolet}{rgb}{0.62, 0.0, 1.0}

\def\rsq{\hspace*{\fill}$\Box$\medskip}

\def\Z{\mathbb Z}
\newtheoremstyle{de}
  {10pt}          
  {10pt}  
  {\rm}  
  {}
  {\bf}  
  {. }    
  { }    
  {}     
\theoremstyle{de}

\newtheorem{example}{Example}[section]

\newtheorem{problem}{Problem}[section]

\newtheoremstyle{theorem}
  {10pt}          
  {10pt}  
  {\it}  
  {}
  {\bf}  
  {. }    
  { }    
  {}     
\theoremstyle{theorem}

\newtheorem{theorem}{Theorem}[section]
\newtheorem{lemma}[theorem]{Lemma}

\newtheorem{corollary}[theorem]{Corollary}

\newtheorem{conjecture}{Conjecture}[section]

\numberwithin{equation}{section}
\def\Z{\mathbb{Z}}

\usepackage[mathlines]{lineno}
\modulolinenumbers[2]
\newcommand*\patchAmsMathEnvironmentForLineno[1]{%
\expandafter\let\csname old#1\expandafter\endcsname\csname #1\endcsname  \expandafter\let\csname oldend#1\expandafter\endcsname\csname end#1\endcsname  \renewenvironment{#1}%
{\linenomath\csname old#1\endcsname}%
{\csname oldend#1\endcsname\endlinenomath}}%
\newcommand*\patchBothAmsMathEnvironmentsForLineno[1]{%
\patchAmsMathEnvironmentForLineno{#1}%
\patchAmsMathEnvironmentForLineno{#1*}}%
\AtBeginDocument{%
\patchBothAmsMathEnvironmentsForLineno{equation}%
\patchBothAmsMathEnvironmentsForLineno{align}%
\patchBothAmsMathEnvironmentsForLineno{flalign}%
\patchBothAmsMathEnvironmentsForLineno{alignat}%
\patchBothAmsMathEnvironmentsForLineno{gather}%
\patchBothAmsMathEnvironmentsForLineno{multline}%
}

\begin{document}
\begin{center}
{\mathversion{bold}\Large \bf On  local antimagic chromatic number of the join of two special families of graphs | II}

\bigskip
{\large  Gee-Choon Lau$^a$, Wai Chee Shiu$^b$ }\\

\medskip

\emph{{$^a$}77D, Jalan Suboh, 85000 Segamat, Johor, Malaysia}\\
\emph{geeclau@yahoo.com}\\

\medskip
\emph{{$^b$}Department of Mathematics,}\\
\emph{The Chinese University of Hong Kong,}\\
\emph{Shatin, Hong Kong, P.R. China.}\\
\emph{wcshiu@associate.hkbu.edu.hk}\\

\end{center}

\begin{abstract}
 It is known that null graphs and 1-regular graphs are the only regular graphs without local antimagic chromatic number.  In this paper, we proved that the join of 1-regular graph and a null graph has local antimagic chromatic number is 3. Consequently, we also obtained many families of (possibly disconnected or regular) bipartite and tripartite graph with local antimagic chromatic number 3. 
\ms


\end{abstract}

\baselineskip18truept
\normalsize

\section{Introduction}
Let $G=(V, E)$ be a connected graph of order $p$ and size $q$.
A bijection $f: E\rightarrow \{1, 2, \dots, q\}$ is called a \textit{local antimagic labeling}
if $f^{+}(u)\neq f^{+}(v)$ whenever $uv\in E$,
where $f^{+}(u)=\sum_{e\in E(u)}f(e)$ and $E(u)$ is the set of edges incident to $u$.
The mapping $f^{+}$ which is also denoted by $f^+_G$ is called a \textit{vertex labeling of $G$ induced by $f$}, and the labels assigned to vertices are called \textit{induced colors} under $f$.
The \textit{color number} of a local antimagic labeling $f$ is the number of distinct induced colors under $f$, denoted by $c(f)$.  Moreover, $f$ is called a \textit{local antimagic $c(f)$-coloring} and $G$ is {\it local antimagic $c(f)$-colorable}. The \textit{local antimagic chromatic number} $\chi_{la}(G)$ is defined to be the minimum number of colors taken over all colorings of $G$ induced by local antimagic labelings of $G$~\cite{Arumugam}. Let $G+H$ and $mG$ denote the disjoint union of graphs $G$ and $H$, and $m$ copies of $G$, respectively. For integers $c < d$, let $[c,d] = \{n\in\Z\;|\; c\le n\le d\}$. Very few results on the local antimagic chromatic number of regular graphs are known (see~\cite{Arumugam, LauLiShiu}).
Throughout this paper, we let $V(aP_2\vee O_m) = \{u_i, v_i, x_j\;|\; 1\le i\le a, 1\le j\le m\}$ and $E(aP_2 \vee O_m) = \{u_ix_j, v_ix_j, u_iv_i\;|\; 1\le i\le a, 1\le j\le m\}$. We also let $V(a(P_2\vee O_m)) = \{u_i, v_i, x_{i,j}\;|\; 1\le i\le a, 1\le j\le m\}$ and $E(a(P_2\vee O_m)) = \{u_ix_{i,j}, v_ix_{i,j}, u_iv_i\;|\; 1\le i\le a, 1\le j\le m\}$.

\ms \nt In~\cite{Haslegrave}, the author proved that all connected graphs without a $P_2$ component admit a local antimagic labeling. Thus, $O_m, m\ge 1$ and $aP_2, a\ge 1$ are the only families of regular graphs without local antimagic chromatic number.  In~\cite{Arumugam}, it was shown that $\chi_{la}(aP_2\vee O_1) = 3$ for $a\ge 1$. In~\cite{LauShiu-join, LauShiu-chila3}, the authors proved that $\chi_{la}((2k+1)P_2 \vee O_m)=3$ for all $k\ge 1, m\ge 2$.  In this paper, we extend the ideas in~\cite{LauShiu-join, LauShiu-chila3, LSNP-even, LSPN-odd} to further  prove that $\chi_{la}((2k)P_2 \vee O_m) = 3$ for all $k\ge 1$, $m\ge 2$. Moreover, we also obtain other families of  bipartite and tripartite graphs with local antimagic chromatic number 3.

\ms\nt The following lemma in \cite[Lemma 2.1]{LSN-DMGT} or~\cite[Lemma 2.3]{LSN-IJMSI} is needed.

\begin{lemma}[{\cite[Lemma 2.3]{LSN-IJMSI}}]\label{lem-2part} Let $G$ be a graph of size $q$. Suppose there is a local antimagic labeling of $G$ inducing a $2$-coloring of $G$ with colors $x$ and $y$, where $x<y$. Let $X$ and $Y$ be the sets of vertices colored $x$ and $y$, respectively, then $G$ is a bipartite graph with bipartition $(X,Y)$ and $|X|>|Y|$. Moreover,
$x|X|=y|Y|= \frac{q(q+1)}{2}$.
\end{lemma}

\nt The contrapositive implies that a connected bipartite graph $G$ with equal partite set size must have $\chi_{la}(G)\ge 3$. Note that if $G$ is a disconnected bipartite graph with each component having equal partite set size, then we also have $\chi_{la}(G) \ge 3$.

\section{Matrix of size $(4n+1)\times 2k$} \label{sec:2k}

Consider $(2k)(P_2 \vee O_{2n})$ of order $2k(2n+2)$ and size $2k(4n+1)$ for $k, n\ge 1$. We shall construct the following tables, which shows the label of each edge under a labeling $f$. First, we assume $k\ge 2$.

\begin{table}[H]
\fontsize{7}{10}\selectfont
$\begin{tabu}{|c|[1pt]c|[1pt]c|c|c|c|c|[1pt]}\hline
i & 1 & 2 & 3 & \cdots & k-1 & k \\\tabucline[1pt]{-}
f(u_ix_{i,2n-1})  & 4k(n-1) + 8k+1 & 4k(n-1) + 9k+1 & 4k(n-1) + 9k+2 &  \cdots & 4k(n-1) + 10k-2 & 4k(n-1) + 10k-1 \\\hline
f(u_ix_{i,2n})   & 4k(n-1) + 6k & 4k(n-1) + 6k+1 & 4k(n-1) + 6k+2 & \cdots &  4k(n-1) + 7k-2 & 4k(n-1) + 7k-1 \\\tabucline[1pt]{-}
f(u_iv_i) & 4k(n-1) + 7k & 4k(n-1) + 6k-1 & 4k(n-1) + 6k-3 & \cdots  & 4k(n-1) + 4k+5  & 4k(n-1) + 4k+3 \\\tabucline[1pt]{-}
\end{tabu}$\\[1mm]

$\begin{tabu}{|c|[1pt]c|c|c|c|c|[1pt]c|}\hline
i & k+1 & k+2 & \cdots & 2k-2 &  2k-1 & 2k \\\tabucline[1pt]{-}
f(u_ix_{i,2n-1})  & 4k(n-1) + 8k+2 & 4k(n-1) + 8k+3 &  \cdots & 4k(n-1) + 9k-1 & 4k(n-1) + 9k & 4k(n-1)+10k \\\hline
f(u_ix_{i,2n})   & 4k(n-1) + 7k+1 & 4k(n-1) + 7k+2 & \cdots &  4k(n-1) + 8k-2 & 4k(n-1) + 8k-1 & 4k(n-1) + 8k \\\tabucline[1pt]{-}
f(u_iv_i) & 4k(n-1) + 6k-2 & 4k(n-1) + 6k-4 & \cdots  & 4k(n-1) + 4k+4  & 4k(n-1) + 4k+2 & 4k(n-1) + 3k+1 \\\tabucline[1pt]{-}
\end{tabu}$
\caption{Each column sum is $12kn+9k+1$.}\label{T-1}
\end{table}

\begin{table}[H]
\fontsize{7}{10}\selectfont
$\begin{tabu}{|c|[1pt]c|[1pt]c|c|c|c|c|[1pt]}\hline
i & 1 & 2 & 3 & \cdots & k-1 & k \\\tabucline[1pt]{-}
f(u_iv_i) & 4k(n-1) + 7k & 4k(n-1) + 6k-1 & 4k(n-1) + 6k-3 & \cdots  & 4k(n-1) + 4k+5  & 4k(n-1) + 4k+3 \\\tabucline[1pt]{-}
f(v_ix_{i,2n-1})  & 4k(n-1) + 1  & 4k(n-1) + k+1 & 4k(n-1) + k+2 & \cdots  & 4k(n-1) + 2k-2 & 4k(n-1) + 2k-1  \\\hline
f(v_ix_{i,2n})   & 4k(n-1)+2k+1 & 4k(n-1) + 2k+2 & 4k(n-1) + 2k+3 & \cdots   & 4k(n-1)+ 3k-1 &  4k(n-1) + 3k \\\hline
\end{tabu}$\\[1mm]

$\begin{tabu}{|c|[1pt]c|c|c|c|c|[1pt]c|}\hline
i & k+1 & k+2 & \cdots & 2k-2 &  2k-1 & 2k \\\tabucline[1pt]{-}
f(u_iv_i) & 4k(n-1) + 6k-2 & 4k(n-1) + 6k-4 & \cdots  & 4k(n-1) + 4k+4  & 4k(n-1) + 4k+2 & 4k(n-1) + 3k+1 \\\tabucline[1pt]{-}
f(v_ix_{i,2n-1})  & 4k(n-1) + 2  & 4k(n-1) + 3 &  \cdots  & 4k(n-1) + k-1 & 4k(n-1) + k &  4k(n-1) + 2k \\\hline
f(v_ix_{i,2n})   & 4k(n-1)+3k+2 & 4k(n-1) + 3k+3 &  \cdots   & 4k(n-1)+ 4k-1 &  4k(n-1) + 4k & 4k(n-1) + 4k+1 \\\hline
\end{tabu}$
\caption{Each column sum is $12kn-3k+2$.}\label{T-2}
\end{table}

\nt For $j\in[1, n-1]$,
\begin{table}[H]
$\fontsize{7}{10}\selectfont
\begin{tabu}{|c|[1pt]c|c|c|c|c|}\hline
i & 1 & 2  & \cdots  & 2k-1 & 2k \\\tabucline[1pt]{-}
f(u_ix_{i,2j-1}) & 2k(4n+3-2j) -2k+1 & 2k(4n+3-2j) -2k+2  & \cdots  & 2k(4n+3-2j)-1  &   2k(4n+3-2j) \\\hline
f(u_ix_{i,2j}) &   2k(2j-1)+2k &  2k(2j-1)+2k-1  & \cdots  & 2k(2j-1)+2 & 2k(2j-1)+1  \\\hline
\end{tabu}$
\caption{Each column sum is $8kn+4k+1$.}\label{T-3}
\end{table}
\begin{table}[H]
$\fontsize{7}{10}\selectfont
\begin{tabu}{|c|[1pt]c|c|c|c|c|}\hline
i & 1 & 2  & \cdots  & 2k-1 & 2k \\\tabucline[1pt]{-}
f(v_ix_{i,2j-1}) & 4k(j-1) + 1 & 4k(j-1) + 2  & \cdots & 4k(j-1) + 2k-1 & 4k(j-1) + 2k \\\hline
f(v_ix_{i,2j}) & 2k(4n+2-2j) & 2k(4n+2-2j) -1 & \cdots & 2k(4n+2-2j) -2k+2   & 2k(4n+2-2j) -2k+1   \\\hline
\end{tabu}$
\caption{Each column sum is $8kn+1$.}\label{T-4}
\end{table}

\nt When $n=1$, we only use Tables~\ref{T-1} and~\ref{T-2}. For $n\ge 2$ and $k=1$, we only use the first and the last columns of each table. For clarity, the table is given below, where $j\in[1, n-1]$.

\[\fontsize{7}{10}\selectfont
\begin{tabu}{|c|[1pt]c|c|}\hline
i & 1 & 2    \\\tabucline[1pt]{-}
f(u_ix_{i,2j-1})  & 8n+5-4j   & 8n+6-4j   \\\hline
f(u_ix_{i,2j}) &  4j & 4j-1 \\\hline
f(u_ix_{i,2n-1}) & 4n+5 & 4n+6 \\\hline
f(u_ix_{i,2n})  &  4n+2 & 4n+4  \\\tabucline[1pt]{-}
f(u_iv_i) &  4n+3  & 4n  \\\tabucline[1pt]{-}
f(v_ix_{i,2j-1}) &  4j-3 & 4j-2 \\\hline
f(v_ix_{i,2j}) & 8n+4-4j & 8n+3-4j  \\\hline
f(v_ix_{i,2n-1})  &  4n-3 & 4n-2  \\\hline
f(v_ix_{i,2n})   &  4n-1 & 4n+1 \\\hline
\end{tabu}\]

We now have the following observations that hold for $(n,k)\ne (1,1)$.
\begin{enumerate}[(A)]
\item Each integer in $[1, 2k(4n+1)]$ serves as an edge label once.
\item From the above tables, we have $f^+(u_i)=(n-1)(8kn+4k+1)+12kn+9k+1=8kn^2+8kn+5k+n$ and $f^+(v_i)= (n-1)(8kn+1)+12kn-3k+2=8kn^2+4kn-3k+n+1$ for each $i$. Clearly, $f^+(u_i)>f^+(v_i)$.
\item Suppose $n=1$. For each $1\le i\le k$, $f(u_ix_{i,2n-1}) + f(u_{2k+1-i}x_{2k+1-i,2n-1}) = 8kn+10k+1$, $f(u_ix_{i,2n}) + f(u_{2k+1-i}x_{2k+1-i,2n}) = 8kn+6k$, whereas $f(v_ix_{i,2n-1}) + f(v_{2k+1-i}x_{2k+1-i,2n-1}) = 8kn - 6k+1$, $f(v_ix_{i,2n}) + f(v_{2k+1-i}x_{2k+1-i,2n}) = 8kn - 2k+2$. Thus, sum of rows $f(u_ix_{i,2n-1})$ and $f(v_ix_{i,2n-1})$ entries is $k(8kn+10k+1)+k(8kn - 6k+1) = k(16kn+4k+2)$. Similarly, sum of rows $f(u_ix_{i,2n})$ and $f(v_ix_{i,2n})$ entries is $k(8kn+6k) + k(8kn-2k+2) = k(16kn+4k+2)$.
\item Suppose $n\ge 2$ and $1\le j\le n-1$. For each $1\le i\le k$, $f(u_ix_{i,2j-1}) + f(u_{2k+1-i}x_{2k+1-i,2j-1}) =4k(4n+3-2j)-2k+1$, $f(u_ix_{i,2j}) + f(u_{2k+1-i}x_{2k+1-i,2j}) = 4k(2j-1)+2k+1$, whereas $f(v_ix_{i,2j-1}) + f(v_{2k+1-i}x_{2k+1-i,2j-1}) = 8k(j-1)+2k+1$, $f(v_ix_{i,2j}) + f(v_{2k+1-i}x_{2k+1-i,2j}) =4k(4n+2-2j) - 2k+1$. Together with Observation (C), we conclude that for $1\le j\le n$, sum of rows $f(u_ix_{i,2j-1})$ and $f(v_ix_{i,2j-1})$ entries is $k[4k(4n+3-2j)-2k+1]+k[8k(j-1)+2k+1] = k(16kn+4k+2)$, and sum of rows $f(u_ix_{i,2j})$ and $f(v_ix_{i,2j})$ entries is $k[4k(2j-1)+2k+1] + k[4k(4n+2-2j) - 2k+1] = k(16kn + 4k + 2)$.
\item Suppose $k = rs$, $r\ge 2$, $s\ge 1$. We can divide the $2k$ columns into $2r$ blocks of $s$ column(s) with the $b$-th block containing $(b-1)s+1, (b-1)s+2, \ldots, bs$ columns. From Observations (C) and (D), we can conclude that sum of row $f(u_ix_{i,j})$ and $f(v_ix_{i,j})$ entries in block $b$-th and $(2r+1-b)$-th collectively is a constant $s(16kn+4k+2)$ for $1\le j\le 2n, 1\le b\le r, s\ge 1$.
\item For each $1\le i\le k$ and $1\le j\le 2n$, $f(u_ix_{i,j}) + f(v_{2k+1-i}x_{2k+1-i,j}) = f(v_ix_{i,j}) + f(u_{2k+1-i}x_{2k+1-i,j}) = 8kn+2k+1$.
\end{enumerate}

\begin{theorem}\label{thm-evenP2VO2n} For $k, n\ge 1$, $\chi_{la}((2k)P_2 \vee O_{2n}) = 3$. \end{theorem}

\begin{proof} Note that $\chi_{la}((2k)P_2 \vee O_{2n}) \ge \chi((2k)P_2\vee O_{2n}) = 3$. A local antimagic 3-coloring of $2P_2\vee O_2$ is given in the figure below with induced vertex labels $14, 19, 22$.

\begin{figure}[H]
\centerline{\epsfig{file=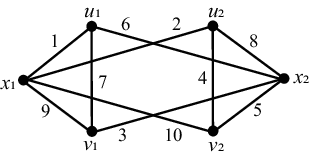, width=5cm}}
\caption{Graph $2P_2\vee O_2$.}\label{fig:2P2VO2}
\end{figure}

\ms\nt Suppose $(n,k)\ne (1,1)$. Consider $G = 2k(P_2\vee O_{2n})$. Since $G$ has size $2k(4n+1)$, we can now define a bijection $f : E(G) \to [1, 2k(4n+1)]$ according to the tables above. For each $j\in [1,2n]$, merging the vertices in $\{x_{i,j} \;|\; 1\le i\le 2k\}$ to form new vertex $x_j$ of degree $4k$, we get the graph $(2k)P_2\vee O_{2n}$. We still use $f$ to denote the labeling for the graph $(2k)P_2\vee O_{2n}$. From Observations (A) to (D) above, we get that $(2k)P_2 \vee O_{2n}$ admits a bijective edge labeling $f$ with
\begin{multicols}{2}
\begin{enumerate}[$(a)$]
\item $f^+(x_j) = k(16kn+4k+2)$,
\item $f^+(u_i) = 8kn^2+8kn+5k+n$,
\item $f^+(v_i) =8kn^2+4kn-3k+n+1$,
\item $f^+(u_i)>f^+(v_i)$.
\end{enumerate}
\end{multicols}

\nt Now,
\begin{eqnarray*}
(a)-(b) &=& 16k^2n - 8kn^2 + 4k^2 - 6kn - 3k - n \\
&=& (8kn+2k)(2k-n) - 6kn - 3k -n \\
&<& 0 \quad \mbox{ if $2k\le n$}.
\end{eqnarray*}
Otherwise, $2k\ge n+1$ so that $(a) - (b)\ge 2kn - k - n > 0$ since $(n,k)\ne (1,1)$. Thus, $(a) \ne (b)$. Similarly,
\begin{eqnarray*}
(a)-(c) &=& 16k^2n-4kn^2+4k^2-8kn+5k-n-1\\
&=&(8kn + 2k)(2k-n) - k(6n +5)-n-1 \\
&<& 0 \quad \mbox{ if $2k\le n$}.
\end{eqnarray*}
If $2k = n+1$, then $(a) - (c) = 2kn - 3k - n -1 = 2k(2k-1)-3k-(2k-1)-1 = 4k^2 - 7k\ne 0$ Otherwise, $2k\ge n+2$ so that $a) - (c) \ge 10kn - k - n -1 > 0$. Thus, $(a) \ne (c)$.

\ms\nt Therefore, $f$ is a local antimagic 3-coloring and $\chi_{la}((2k)P_2\vee O_{2n})\le 3$. This completes the proof.
\end{proof}

\begin{example} For $n= 2, k =4$ so that $j=1$, we have the following $9 \times 8$ matrix and $8(P_2 \vee O_4)$ as shown. For each $j\in[1,4]$, merge the vertices in $\{x_{i,j}\mid 1\le i\le 8\}$, we get the $8P_2\vee O_4$ and the labeling as defined above.

\[\fontsize{7}{10}\selectfont
\begin{tabu}{|c|[1pt]c|c|c|c|c|c|c|c|}\hline
i & 1 & 2 & 3 & 4 & 5 & 6 & 7 & 8  \\\tabucline[1pt]{-}
f(u_ix_{i,1}) &65 & 66 & 67 & 68 & 69 & 70 & 71 & 72  \\\hline
f(u_ix_{i,2}) & 16 & 15 & 14 & 13 & 12 & 11 & 10 & 9   \\\hline
f(u_ix_{i,3})  & 49 & 53 & 54 & 55 & 50 & 51 &52 & 56  \\\hline
f(u_ix_{i,4})   & 40 & 41 & 42 & 43& 45 & 46 & 47 & 48  \\\tabucline[1pt]{-}
f(u_iv_i) & 44 & 39 & 37 & 35 & 38 & 36 & 34 & 29  \\\tabucline[1pt]{-}
f(v_ix_{i,1}) &  1 & 2 &  3 & 4 & 5 & 6 & 7 & 8  \\\hline
f(v_ix_{i,2}) & 64 & 63 & 62 & 61 & 60 & 59 & 58 & 57    \\\hline
f(v_ix_{i,3})  & 17 & 21 & 22 & 23 & 18 & 19 & 20 & 24    \\\hline
f(v_ix_{i,4})   &  25 & 26 & 27 & 28 & 30 & 31 & 32 & 33   \\\hline
\end{tabu}\]
\begin{figure}[H]
\centerline{\epsfig{file=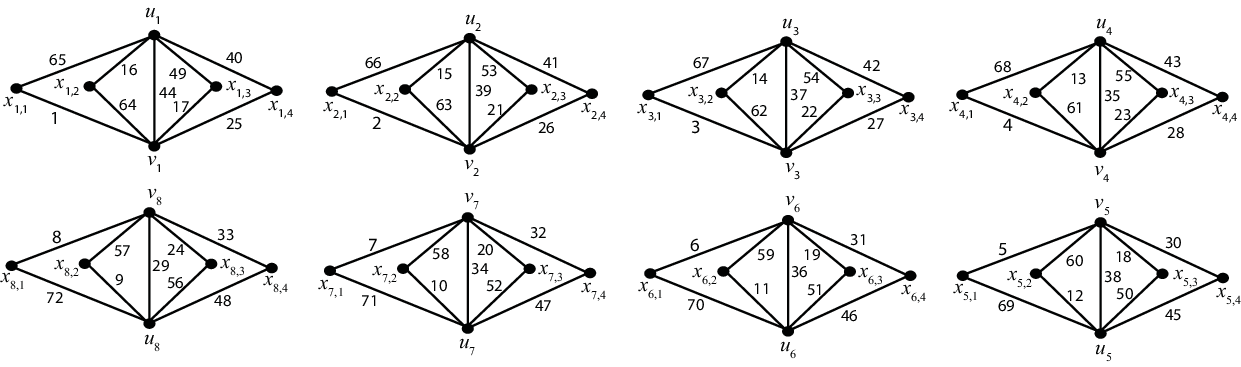, width=15cm}}
\caption{Graph $8(P_2\vee O_4)$.}\label{fig:8(P2VO4))}
\end{figure}
\rsq
\end{example}

\begin{theorem}\label{thm-r(sP2VO2n)} For $n\ge 1, r\ge 2, s\ge 1$, $\chi_{la}(r((2s)P_2\vee O_{2n})) = 3$.  \end{theorem}

\begin{proof} Let $k=rs \ge 2$. Begin with $2k(P_2\vee O_{2n})$ and the local antimagic 3-coloring $f$ as defined in the proof of Theorem~\ref{thm-evenP2VO2n}. By Observations (C) and (D), we can now partition the $2k$ components into $2r$ blocks of $s$ component(s) with the $b$-th block $(1\le b\le 2r)$ containing vertices $u_i, v_i, x_{i,j}$ for $(b-1)s+1 \le i\le bs$, $r\ge 2, s\ge 1$.  Now, for each $j \in [1,2n]$, merge the vertices in $\{x_{i,j}, x_{2k+1-i,j}\;|\; (b-1)s+1\le i \le bs\}$ to form a new vertex of degree $4s$, still denoted $x_{i,j}$,  $1\le i \le r$. We now get an $r$ components graph $r((2s)P_2 \vee O_{2n})$ with a bijective edge labeling, still denoted $f$, such that the degree $2n+1$ vertices, $u_i$ and $v_i$, have induced vertex labels $(a)  = f^+(u_i) =8kn^2+8kn+5k+n$ and $(b) = f^+(v_i) =8kn^2+4kn-3k+n+1$, respectively. Moreover, the vertices $x_{i,j}$ have induced vertex label $(c) = f^+(x_{i,j}) = s(16kn+4k+2)$. Since $(a) > (b)$, we shall show that $(a), (b) \ne (c)$. If $n=1$, the three induced vertex labels are $21k+1, 9k+2, s(20k+2)$ which are distinct. We now assume $n\ge 2$.

\ms\nt Now, $(a)-(c)$ is:
\begin{align}
8kn^2 - 16kns + 8kn - 4ks + 5k + n - 2s = (8kn+1)(n-2s+1) - k(4s-5)-1\label{eq:a-c}
\end{align}

\nt If $n-2s+1\le 0$, then $2\le n < 2s$. Thus $s\ge 2$. We get $\eqref{eq:a-c} < 0$.

\nt If $n-2s+1 >0$, then $n \ge 2s$. Thus,
\begin{align*}
(8kn+1)(n-2s+1)-k(4s-5) - 1 & \ge 8kn +1- 4ks + 5k -1\\
 & = 4k(2n - s) + 5k > 0.
\end{align*}

\nt Next, $(b)-(c)$ is:
\begin{align}
8kn^2 -16kns+4kn - 4ks-3k + n - 2s+1= (4kn+\frac{1}{2})(2n-4s+1)-k(4s+3)+\frac{1}{2}\label{eq:b-c}
\end{align}

\nt If $2n-4s+1\le 0$,  then clearly $\eqref{eq:b-c}<0$.

\nt If $2n-4s+1>0$,  then $n\ge 2s$. Thus,
\begin{align*}
(4kn+\frac{1}{2})(2n-4s+1)-k(4s+3)+\frac{1}{2}\ge 4k(n-s)-3k+1
\ge 4ks-3k+1>0.
\end{align*}

\nt Therefore, $f$ is a local antimagic 3-coloring and $\chi_{la}(r((2s)P_2\vee O_{2n}))\le 3$. This completes the proof.
\end{proof}

\begin{example} Take $n=k=2$ so that $r=2,s=1$, we can get the following $9\times 4$ matrix and $2(2P_2 \vee O_4)$ with the defined labeling as follow. The induced vertex labels are $108, 77, 74$ respectively. 

\[\fontsize{7}{10}\selectfont
\begin{tabu}{|c|[1pt]c|c|c|c|}\hline
i & 1 & 2 & 3 & 4   \\\tabucline[1pt]{-}
f(u_ix_{i,1}) & 33 & 34  & 35  & 36  \\\hline
f(u_ix_{i,2}) &  8  & 7  & 6  & 5   \\\hline
f(u_ix_{i,3})  & 25 & 27 & 26 & 28   \\\hline
f(u_ix_{i,4})   & 20 & 21 & 23 & 24 \\\tabucline[1pt]{-}
f(u_iv_i) & 22 & 19 & 18 & 15 \\\tabucline[1pt]{-}
f(v_ix_{i,1}) &  1 & 2 &  3 & 4  \\\hline
f(v_ix_{i,2}) & 32 & 31 & 30 & 29    \\\hline
f(v_ix_{i,3})  & 9  & 11 & 10 & 12   \\\hline
f(v_ix_{i,4})   & 13 &14 & 16 & 17   \\\hline
\end{tabu}\]
\begin{figure}[H]
\centerline{\epsfig{file=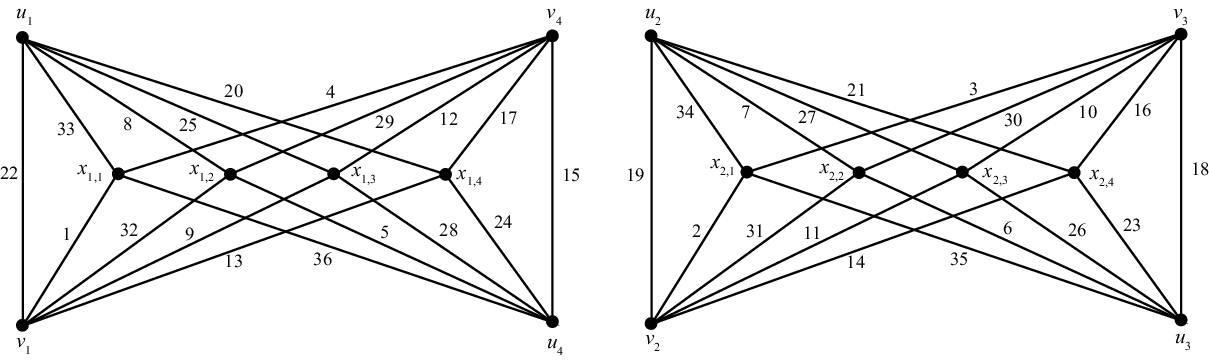, width=14cm}}
\caption{Graph $2(2P_2\vee O_{4})$.}\label{fig:G4(22)}
\end{figure}

\rsq
\end{example}

\nt In the graph $2(2P_2\vee O_4)$ above, we can now split each degree $4s=4$ vertex (with induced vertex label 74) into two degree $2s=2$ vertices with induced vertex label 37. This new graph, say $G_{2,2}(4)$, is bipartite. Each component has a bipartition of equal size $2n+2s=6$. By Lemma~\ref{lem-2part}, $\chi_{la}(G_{2,2}(4))\ge 3$. Since the corresponding edge labeling induces a local antimagic 3-coloring, we have $\chi_{la}(G_{2,2}(4)) = 3$.

\ms\nt Begin with $r((2s)P_2\vee O_{2n})$ and the local antimagic 3-coloring $f$ defined in the proof of Theorem~\ref{thm-r(sP2VO2n)}.  In general, for $k=rs$, $n\ge 1$, $r\ge 2$, $s\ge 1$, let $G_{r,2s}(2n)$ be obtained from $r((2s)P_2\vee O_{2n})$ by splitting each $x_{i,j}$ of degree $4s$ into two vertices, say $y_{i,j}$ and $z_{i,j}$, of degree $2s$ such that both have equal induced vertex labels under $f$. By Observation (F), we know that $G_{r,2s}(2n)$ admits a bijective edge labeling, say $g$, such that $g^+(u_i) = f^+(u_i)=8kn^2+8kn+5k+n$, $g^+(v_i) = f^+(v_i)=8kn^2+4kn-3k+n+1$, and $g^+(y_{i,j}) = g^+(z_{i,j}) = \frac{1}{2}f^+(x_{i,j}) = s(8kn+2k+1)$.

\begin{theorem}\label{thm-G_{r,2s}(2n)} For $n\ge 1$, $r\ge 2$, $s\ge 1$, $\chi_{la}(G_{r,2s}(2n))=3$. \end{theorem}

\begin{proof} By definition, we know that $k = rs\ge 2$ and $G_{r,2s}(2n)$ is  an $r$ components bipartite graph with each component of equal partite sets size $2n+2s$. By Lemma~\ref{lem-2part}, $\chi_{la}(G_{r,2s}(2n))\ge 3$. Moreover, $G_{r,2s}(2n)$ admits a bijective edge labeling $g$ with induced vertex labels $(a) = g^+(u_i) = 8kn^2+8kn+5k+n$, $(b) = g^+(v_i) = 8kn^2+4kn-3k+n+1$, and $(c) = g^+(y_{i,j}) = g^+(z_{i,j})  = s(8kn+2k+1)$. Since $(a) > (b)$, we shall show that $(a), (b) \ne (c)$. Suppose $n=1$. Since $k=rs\ge 2$, we have $(a) = 21k+1$, $(b) = 9k+2$ and $(c) = s(10k+1)$ which are distinct. We now assume $n\ge 2$.

\ms\nt Now, $(a)-(c)$ is:
\begin{align}
8kn^2 - 8kns + 8kn - 2ks + 5k + n - s
= (8kn+1)(n-s+1) - k(2s-5) - 1 \label{eqn:(a)-(c)}
\end{align}

\nt If $n-s+1\le 0$, then  $s\ge n+1\ge 3$. We get $\eqref{eqn:(a)-(c)} < 0$.

\nt If $n-s+1 > 0$, then $n\ge s$. Thus,
\begin{align*}
(8kn+1)(n-s+1) - k(2s-5) - 1 &\ge 8kn+1 - 2ks+5k-1\\
& = 2k(4n - s) + 5k > 0.
\end{align*}

\nt Next, $(b) - (c)$ is:
\begin{align}
8kn^2 - 8kns + 4kn - 2ks - 3k +n - s + 1
 = (4kn+\frac{1}{2})(2n-2s+1) - k(2s+3) + \frac{1}{2} \label{eqn:(b)-(c)}
\end{align}

\nt If $2n-2s+1\le 0$, then clearly $\eqref{eqn:(b)-(c)} < 0$.

\nt If $2n-2s+1 = 1$, then $n=s$. Thus, $(4kn+\frac{1}{2})(2n-2s+1) - k(2s+3) + \frac{1}{2} = k(4n-2s-3) + 1 = k(2s - 3) + 1\ne 0$ for all $k,s$.

\nt If $2n-2s+1\ge 2$, then $2n \ge 2s+1$. Thus,
\begin{align*}
(4kn+\frac{1}{2})(2n-2s+1) - k(2s+3) + \frac{1}{2}  & \ge k(8n-2s-3) + \frac{3}{2} >0.
\end{align*}

\nt Therefore, $g$ is a local antimagic 3-coloring and $\chi_{la}(G_{r,2s}(2n)) \le 3$. This completes the proof.                       \end{proof}

\begin{example} Take $n=2$, $k = 3$ so that $r=3$, $s=1$, the graph $G_{3,2}(4)$ with the local antimagic 3-coloring defined under Theorem~\ref{thm-G_{r,2s}(2n)}  is given in Figure~\ref{fig:G32(4)}.

\begin{figure}[H]
\centerline{\epsfig{file=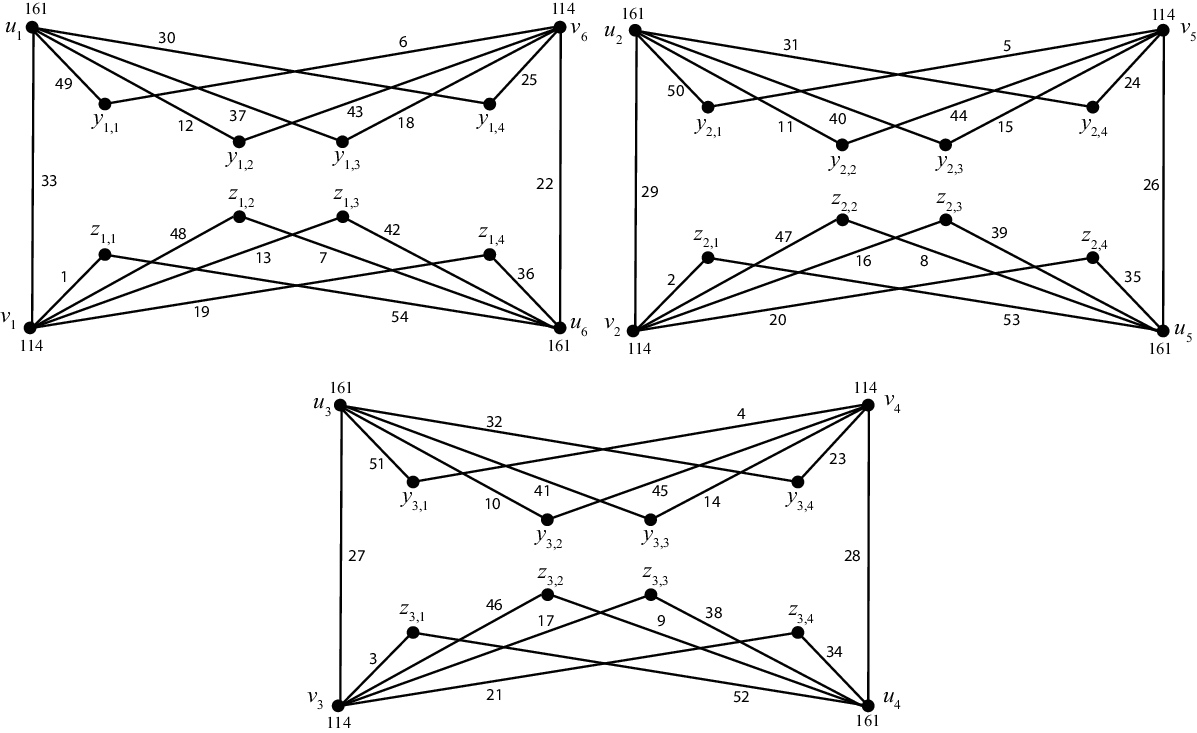, width=13cm}}
\caption{Graph $G_{3,2}(4)$.}\label{fig:G32(4)}
\end{figure}\rsq
\end{example}

\ms\nt Another way to get a new graph is the delete-add process which the authors used in \cite{LauShiu-join, LauShiu-chila3}. In the above example, we can delete edges $v_4x_{1,1}, v_1x_{1,1}$ and $v_3x_{2,1}, v_2x_{2,1}$, and then add edges $v_4x_{2,1}, v_1x_{2,1}$ with labels $4,1$ respectively, also add edges $v_3x_{1,1}, v_2x_{1,1}$ with labels $3,2$ respectively. The new graph is connected but not a graph we have obtained earlier.  In general, let $\mathcal G_{r,2s}(2n)$ be the family of all non-isomorphic graphs that can be obtained by applying the delete-add process to $r((2s)P_2\vee O_{2n})$ and $G_{r,2s}(2n)$. Note that we must have $s\ge 2$ to apply this process to $G_{r,2s}(2n)$. We immediately have the following corollary.

\begin{corollary}\label{cor-G_r,2s(2n)} For $n\ge 1, r\ge 2, s\ge 1$, every graph $G \in \mathcal G_{r,2s}(2n)$ has $\chi_{la}(G) = 3$. \end{corollary}

\nt Recall that in ~\cite[Theorem 2.2]{LauShiu-chila3}, the authors defined a family of $r+1$ componenets tripartite graphs $G_{2n}(2r+1,2s+1)$ obtained from $(2k+1)(P_2 \vee O_{2n})$. Let $2k = (2r+1)(2s)$ for $r, s\ge 1$, we can also define a similar family of $(r+1)$ components tripartite graphs $G_{2n}(2r+1,2s)$ obtained from $2k(P_2 \vee O_{2n})$ by keeping the local antimagic 3-coloring $f$ defined under Theorem~\ref{thm-evenP2VO2n}.

\begin{theorem}\label{thm-G2n(2r+1,2s)} Suppose $n\ge 1$. If $2k = (2r+1)(2s)$ for $r,s\ge 1$, then $\chi_{la}(G_{2n}(2r+1,2s)) = 3$. \end{theorem}

\begin{proof} By definition, $\chi_{la}(G_{2n}(2r+1,2s))\ge \chi(G_{2n}(2r+1,2s)) = 3$. Moroever, we can conclude that $G_{2n}(2r+1,2s)$ admits a bijective edge labeling, say $g$, with induced vertex labels $(a) = 8kn^2+8kn+5k+n$, $(b) = 8kn^2+4kn-3k+n+1$ and $(c) = s(16kn+4k+2)$ such that adjacent vertices must have distinct induced vertex labels. From the proof of Theorem~\ref{thm-r(sP2VO2n)}, we can conclude that $g$ is a local antimagic 3-coloring so that $\chi_{la}(G_{2n}(2r+1,2s))\le 3$. This completes the proof.
\end{proof}

\begin{example} We now use $n=2, k=3$ so that $r= s=1$. We get the $9\times 6$ matrix and the graph $G_4(3,2)$ in Figure~\ref{fig:G4(32)} as shown below.

\[\fontsize{7}{10}\selectfont
\begin{tabu}{|c|[1pt]c|c|c|c|c|c|}\hline
i & 1 & 2 & 3 & 4 & 5 & 6  \\\tabucline[1pt]{-}
f(u_ix_{i,1}) & 49 & 50 & 51 & 52 & 53 & 54  \\\hline
f(u_ix_{i,2}) & 12 & 11 & 10 & 9 &  8  & 7    \\\hline
f(u_ix_{i,3})  & 37 & 40 & 41 & 38 & 39 & 42   \\\hline
f(u_ix_{i,4})   & 30 & 31 & 32 & 34 & 35 & 36  \\\tabucline[1pt]{-}
f(u_iv_i) & 33  & 29 & 27 & 28 & 26 & 22 \\\tabucline[1pt]{-}
f(v_ix_{i,1}) &  1 & 2 &  3 & 4 & 5 & 6  \\\hline
f(v_ix_{i,2}) & 48 & 47 & 46 & 45 & 44 & 43    \\\hline
f(v_ix_{i,3})  & 13  & 16 & 17 & 14 & 15 & 18   \\\hline
f(v_ix_{i,4})   & 19 & 20 & 21 & 23 & 24 & 25   \\\hline
\end{tabu}\]
\begin{figure}[H]
\centerline{\epsfig{file=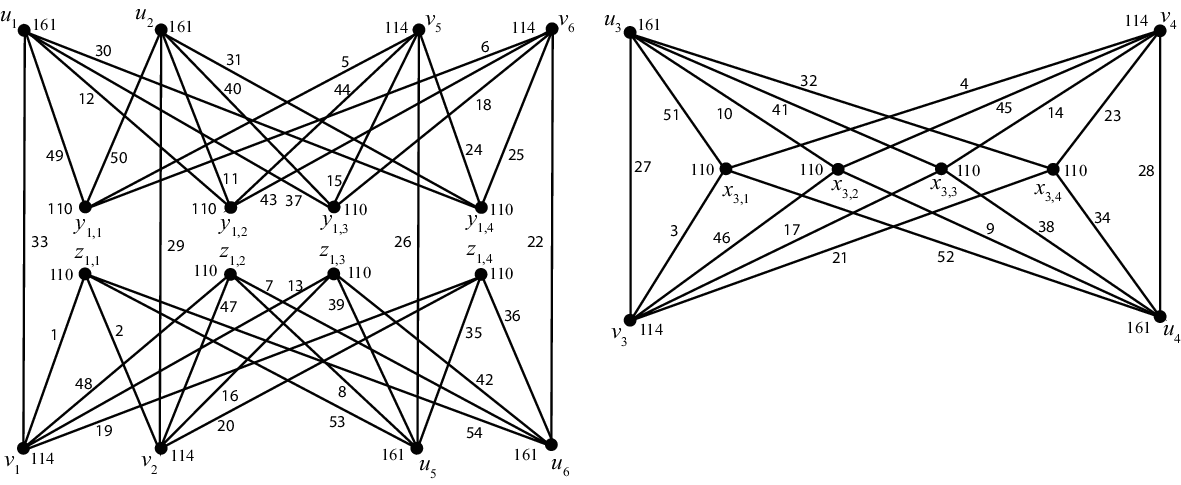, width=13cm}}
\caption{Graph $G_4(3,2)$ with a local antimagic 3-coloring.}\label{fig:G4(32)}
\end{figure}\rsq
\end{example}

\nt We can also apply the delete-add process to $G_{2n}(2r+1,2s)$ to obtain various connected but non-isomorpihc graphs of same order and size, all of which have local antimagic chromatic number 3. In the example above, we may choose to delete edges $u_1y_{1,4}, v_6y_{1,4}, u_3x_{3,1}, v_4x_{3,1}$ and then add edges $u_1x_{3,1}, v_6x_{3,1}, u_3y_{1,4}, v_4y_{1,4}$ with labels $30,25,51,4$ respectively. Let $\mathcal G_{2n}(2r+1,2s)$ be the families of all such graphs. We immediately have the following corollary. \rsq

\begin{corollary}\label{cor-G2n(2r+1,2s)} For $n, r,s\ge 1$, every graph $G\in \mathcal G_{2n}(2r+1,2s)$ has $\chi_{la}(G) = 3$.
\end{corollary}

\nt Note that $k(2P_2\vee O_{2n})$  in Theorem~\ref{thm-r(sP2VO2n)} (respectively $G_{k,2}(2n)$ in Theorem~\ref{thm-G_{r,2s}(2n)}) has $2k$ vertices $v_1, \ldots, v_{2k}$, with induced vertex label $8kn^2+4kn-3k +n +1$. Suppose we can partition $\{v_i\mid 1\le i\le 2k\}$ into blocks of equal size $s\ge 2$ such that all the vertices in the same block do not have common neighbors.
\begin{enumerate}[(a)]
\item Let $\mathcal J_1(k,2n,s)$ be the family of graphs obtained from $k(2P_2\vee O_{2n})$ by merging all the vertices in the same block.
\item Let $\mathcal J_2(k,2n,s)$ be the family of graphs obtained from $G_{k,2}(2n)$ by merging all the vertices in the same block.
\end{enumerate}

\begin{theorem}\label{thm-J1} For $k, s \ge 2$, $n\ge 1$, each graph $G\in \mathcal J_1(k,2n,s)$ has $\chi_{la}(G) = 3$.  \end{theorem}

\begin{proof}  By definition, $\chi_{la}(G) \ge \chi(G) = 3$. Keeping the bijective edge labeling defined for $k(2P_2 \vee O_{2n})$, we can immediately conclude that $G$ admits a bijective labeling with induced vertex labels $(a) = 8kn^2 + 8kn + 5k + n$, $(b) = s(8kn^2 + 4kn - 3k + n + 1)$ and $(c) = 16kn+ 4k + 2$.   It is easy to verify that all the induced vertex labels are distinct and the labeling is local antimagic. Thus, $\chi_{la}(G) \le 3$. Consequently, $\chi_{la}(G) = 3$.  \end{proof}

\nt Similarly, we also have the following theorem with the proof omitted.

\begin{theorem}\label{thm-J2} Suppose $k,s\ge 2$, $n\ge 1$ and $G\in \mathcal J_2(k,2n, s)$. If $G$ is bipartite with every component of equal partite set size or $G$ is tripartite, then $\chi_{la}(G) = 3$. \end{theorem}

\nt The $i$-th component of the graph $k(2P_2\vee O_{2n})$ or graph $G_{k,2}(2n)$ has vertices $u_i$, $v_i$, $u_{2k+1-i}$, $v_{2k+1-i}$, $1\le i\le k$. We now define two new families of connected graphs.
\begin{enumerate}[(a)]
\item Let $J_1(k,2n)$ be a connected graph obtained from $k(2P_2\vee O_{2n})$ by merging vertices in $\{v_i,v_{2k-i}\;|\; 1\le i\le k-1\} \cup \{v_k, v_{2k}\}$. Here $J_1(k,2n)\in \mathcal J_1(k, 2n, 2)$.
\item Let $J_2(k,2n)$ be a connected graph obtained from $G_{k,2}(2n)$ by merging vertices in\break $\{v_i,v_{2k-i}\;|\; 1\le i\le k-1\} \cup \{v_k, v_{2k}\}$. Here $J_2(k,2n)\in \mathcal J_2(k, 2n, 2)$.
\end{enumerate}

\begin{example} Using $6(2P_2\vee O_{2})$, we get partition $\{v_i\;|\; 1\le i\le 12\}$ into 4 blocks of size 3 like $\{v_{3a-2}, v_{3a-1}, v_{3a}\;|\; 1\le a\le 4\}$ to get a 2 components graph, or else $\{v_1, v_2, v_3\}$, $\{v_8, v_9, v_{10}\}$, $\{v_5, v_6, v_7\}$ and $\{v_4, v_{11}, v_{12}\}$ to get a connected graph in $\mathcal J_1(6, 2, 3)$.  Keeping the edge labeling, the induced vertex labels are $127, 168, 122$. \rsq
\end{example}

\ms\nt Suppose we group the $k$ (isomorphic) components of $k(2P_2\vee O_{2n})$ into $t$ graphs, $t\ge 2$. Let the $a$-th graph be $J_a$ which contains $k_a\ge 2$ components. By a similar merging process as getting $J_1(k, 2n)$, we obtain a connected graph, denoted by $H_1(k_a, 2n)$.
So that the new graph $H_1(k_1,2n) + H_1(k_2,2n) + \cdots + H_1(k_t,2n)$ is a $t$ components tripartite graph with $k = k_1 + \cdots + k_t$. Thus, by keeping the edge labeling of $k(2P_2 \vee O_{2n})$, we have the following corollary immediately.

\begin{corollary}\label{cor-H1} For $n\ge 1$ and $t, k_a\ge 2$, the graph $G=H_1(k_1,2n) + H_1(k_2,2n) + \cdots + H_1(k_t,2n)$ has $\chi_{la}(G) = 3$. \end{corollary}

\nt Using $G_{k,2}(2n)$, we can also construct the graph $G = H_2(k_1,2n) + H_2(k_2,2n) + \cdots + H_2(k_t,2n)$ which is tripartite if some of $k_i$ is odd.    Otherwise,  each $H_2(k_i,2n)$ is a bipartite graph with equal partite set size $2k_in + \frac{3k_i}{2}$. Thus, $G$ is a bipartite graph with equal partite set size $\sum\limits^t_{i=1} (2k_in + \frac{3k_i}{2})$. We also have the following corollary immediately.

\begin{corollary}\label{cor-H2} For $n\ge 1$ and $t, k_a\ge 2$, the graph $G=H_2(k_1,2n) + H_2(k_2,2n) + \cdots + H_2(k_t,2n)$ has $\chi_{la}(G) = 3$. \end{corollary}

\begin{example} For example, take $6(2P_2\vee O_{2})$, we can get either $H_1(2,2) + H_1(4,2)$ if we merge the vertices in $\{v_1,v_2\}$, $\{v_{11}, v_{12}\}$, $\{v_3, v_9\}$, $\{v_4, v_8\}$, $\{v_5, v_7\}$ and $\{v_6,v_{10}\}$,  or else we get $2H_1(3,2n)$ if we merge the vertices in $\{v_1, v_{11}\}$, $\{v_2, v_{10}\}$, $\{v_3, v_{12}\}$, $\{v_4, v_8\}$, $\{v_5, v_7\}$ and $\{v_6,v_9\}$. The induced vertex labels are $127, 112, 122$. \rsq
\end{example}

\section{Matrix of size $(4n+3)\times 2k$ }

\nt Consider $(2k)(P_2\vee O_{2n+1})$ of order $2k(2n+3)$ and size $2k(4n+3)$ for $k,n\ge 1$. We shall construct the following table, which shows the label of each edge under a labeling $f$.  
\begin{table}[H]
\fontsize{7}{10}\selectfont
$\begin{tabu}{|c|[1pt]c|c|c|c|c|[1pt]}\hline
i & 1 & 2  & \cdots & 2k-1 & 2k  \\\tabucline[1pt]{-}
f(u_ix_{i,1}) & 4n(2n-1)+10k & 4n(2n-1)+10k-1 & \cdots & 4n(2n-1)+8k+2 & 4n(2n-1)+8k+1\\\hline
f(u_ix_{i,2}) & 4n(2n-1)+6k+1 & 4n(2n-1)+6k+2 & \cdots & 4n(2n-1)+8k-1 & 4n(2n-1)+8k\\\hline
\vdots & \vdots &  \vdots & \cdots & \vdots & \vdots   \\\hline
f(u_ix_{i,2j-1}) & 4k(2n-j)+10k & 4k(2n-j) +10k-1& \cdots &  4k(2n-j)+8k+2 &  4k(2n-j) + 8k+1  \\\hline
f(u_ix_{i,2j}) &  4k(2n-j)  +6k+1& 4k(2n-j) + 6k+2 & \cdots & 4k(2n-j)  +8k-1 &  4k(2n-j) +  8k \\\hline
\vdots & \vdots &  \vdots & \cdots & \vdots & \vdots   \\\hline
f(u_ix_{i,2n-1}) & 4kn+10k & 4kn+10k-1 & \cdots & 4kn+8k+2 & 4kn+8k+1 \\\hline
f(u_ix_{i,2n}) & 4kn+6k+1 & 4kn+6k+2 & \cdots & 4kn+8k-1 & 4kn+8k \\\hline
f(u_ix_{i,2n+1}) & 4kn+6k & 4kn+ 6k-1  & \cdots & 4kn+4k+2 & 4kn+4k+1 \\\tabucline[1pt]{-}
f(u_iv_i) & 1 & 2  & \cdots & 2k-1 & 2k  \\\tabucline[1pt]{-}
f(v_ix_{i,1}) & 4k & 4k-1 & \cdots & 2k+2 & 2k+1    \\\hline
f(v_ix_{i,2}) & 4k+1 & 4k+2 & \cdots & 6k-1 & 6k \\\hline
f(v_ix_{i,3}) & 8k & 8k-1 & \cdots & 6k+2 & 6k+1 \\\hline
\vdots & \vdots & \vdots & \cdots & \vdots & \vdots   \\\hline
f(v_ix_{i,2j}) &  4kj+1  & 4kj+2 & \cdots & 4kj+2k-1  & 4kj + 2k    \\\hline
f(v_ix_{i,2j+1}) & 4kj+4k  & 4kj+4k-1   & \cdots & 4kj+2k+2  & 4kj+2k+1    \\\hline
\vdots & \vdots & \vdots & \cdots & \vdots & \vdots   \\\hline
f(v_ix_{i,2n+1}) & 4kn+4k & 4kn+4k-1 & \cdots & 4kn+2k+2 & 4kn+2k+1 \\\hline
\end{tabu}$
\caption{$(4n+3)\times 2k$ matrix for  $(2k)(P_2 \vee O_{2n+1})$}.
\end{table}

\nt We now have the following observations.
\begin{enumerate}[(1)]
\item Each integer in $[1,2k(4n+3)]$ serves as an edge label once.
\item For each column, the sum of the first $2n+2$ entries is\\  $f^+(u_i)=4kn+6k+1 + \sum\limits^{n}_{j=1} [8k(2n-j)+16k+1] = (n+1)(12nk+4k+1)+2k$.
\item For each column, the sum of the last $2n+2$ entries is\\ $f^+(v_i)=4k+1 + \sum\limits^{n}_{j=1} [4k+1 + 8kj]= (n+1)(4nk+4k+1)$. Clearly $f^+(u_i)>f^+(v_i)$.
\item Consider $1\le i\le 2k$. For odd $j\in [1,2n+1]$, the terms $f(u_ix_{i,j}) + f(v_ix_{i,j})$ form an arithmetic sequence with first term $8kn+10k$, last term $8kn+6k+2$ and common difference $-2$. For even $j\in [1,2n+1]$, the terms $f(u_ix_{i,j}) + f(v_ix_{i,j})$ form an arithmetic sequence with first term $8kn+6k+2$, last term $8kn+10k$ and common difference $2$. Thus, for each $j\in[1,2n+1]$, sum of rows $f(u_ix_{i,j})$ and $f(v_ix_{i,j})$ entries is $2k(8kn+8k+1)$.

 Note that, if $\{A_i\}_{i=1}^m$ is an increasing (resp. decreasing) sequence, then $\{A_{m+1-i}\}_{i=1}^m$ is a decreasing (resp. increasing) sequence.

    Thus, for $1\le i\le k$, $f(u_ix_{i,j}) + f(v_ix_{i,j}) + f(u_{2k+1-i}x_{2k+1-i,j}) + f(v_{2k+1-i}x_{2k+1-i,j}) = K$, where $K$ is independent of $i$. By the discussion above, $K=    2(8kn+8k+1)$.

\item Suppose $k=rs$ for $r\ge 2, s\ge 1$. We can divide the $2k$ columns into $2r$ blocks of $s$ column(s) with the $b$-th block containing $(b-1)s+1, (b-1)s+2, \ldots, bs$ columns. From Observation~(3), we can conclude that sum of rows $f(u_ix_{i,j})$ and $f(v_ix_{i,j})$ entries in block $b$-th and $(2r+1-b)$-th collectively is a constant $s(16kn+16k+2)$.
\item For each $1\le i\le k$ and $1\le j\le 2n+1$, $f(u_ix_{i,j}) + f(v_{2k+1-i}x_{2k+1-i,j}) = f(v_ix_{i,j}) + f(u_{2k+1-i}x_{2k+1-i,j}) = 8kn+8k+1$.

\end{enumerate}

\begin{theorem}\label{thm-evenP2VO2n+1} For $n,k\ge 1$, $\chi_{la}((2k)P_2 \vee O_{2n+1}) = 3$. \end{theorem}

\begin{proof} Note that $\chi_{la}((2k)P_2\vee O_{2n+1}) \ge \chi((2k)P_2\vee O_{2n+1}) = 3$. Consider $(2k)(P_2 \vee O_{2n+1})$.  Since $G$ has size $2k(4n+3)$, we can now define a bijection $f : E(G) \to [1, 2k(4n+3)]$ according to the table above. Clearly, for $1\le i\le 2k$, $f^+(u_i) = (n+1)(12nk+4k+1)+2k > f^+(v_i) = (n+1)(4nk+4k+1)$. Now, for each $j\in [1,2n+1]$, merging the vertices in $\{x_{i,j} \;|\; 1\le i\le 2k\}$, to form new vertex $x_j$ of degree $4k$, we get the graph $(2k)P_2 \vee O_{2n+1}$. From Observations (1) to (3) above, we get that $(2k)P_2 \vee O_{2n+1}$ admits a bijective edge labeling with
\begin{enumerate}[$(a)$]
\item $f^+(x_j) = 2k(8nk+8k+1) = 16k^2(n+1) + 2k$,
\item $f^+(u_i) = (n+1)(12nk+4k+1)+2k$, and
\item $f^+(v_i) = (n+1)(4nk+4k+1)$.
\end{enumerate}

\nt Now, $(a) - (b) = -(n+1)[4k(3n-4k-1)+1] < 0$ if $3n-4k-1\ge 0$. Otherwise, $(a) - (b) > 0$. Thus, $(a) \ne (c)$. Similarly, $(a) - (c) = (4kn+4k)(4k-n)-4kn-2k-n-1< 0$ if $4k-n\le 1$. Otherwise, $(a)-(c) > 0$.  Thus, $(a) \ne (c)$. Therefore, $f$ is a local antimagic 3-coloring and $\chi_{la}((2k)P_2 \vee O_{2n+1})\le 3$.  This completes the proof.     \end{proof}

\begin{example} Take $n=2, k=3$, we have the following $11\times 6$ matrix and $6(P_2 \vee O_{5}$ as shown. For each $j\in [1,5]$, merge the vertices in $\{x_{i,j}\mid 1\le i\le 6\}$, we get a $6P_2 \vee O_5$ and the labeling as defined above.

\[\fontsize{7}{10}\selectfont
\begin{tabu}{|c|[1pt]c|c|c|c|c|c|}\hline
i & 1 & 2 & 3 & 4 & 5 & 6   \\\tabucline[1pt]{-}
f(u_ix_{i,1}) &  66 & 65 & 64 & 63 & 62 & 61   \\\hline
f(u_ix_{i,2}) & 55 & 56 & 57 & 58 & 59 & 60     \\\hline
f(u_ix_{i,3})  & 54 & 53 & 52 & 51 & 50 & 49     \\\hline
f(u_ix_{i,4}) & 43 & 44 & 45 & 46 & 47 & 48     \\\hline
f(u_ix_{i,5}) & 42 & 41 & 40  & 39 & 38 & 37 \\\tabucline[1pt]{-}
f(u_iv_i) & 1 & 2 & 3 & 4 & 5 & 6\\\tabucline[1pt]{-}
f(v_ix_{i,1}) & 12 & 11 & 10 & 9 &  8 & 7  \\\hline
f(v_ix_{i,2}) & 13  & 14 & 15 & 16 & 17 & 18      \\\hline
f(v_ix_{i,3}) & 24 & 23 & 22 & 21 & 20 & 19     \\\hline
f(v_ix_{i,4})  & 25 & 26 & 27 & 28 & 29 & 30    \\\hline
f(v_ix_{i,5}) & 36 & 35 & 34 & 33 & 32 & 31  \\\hline
\end{tabu}\]

\begin{figure}[H]
\centerline{\epsfig{file=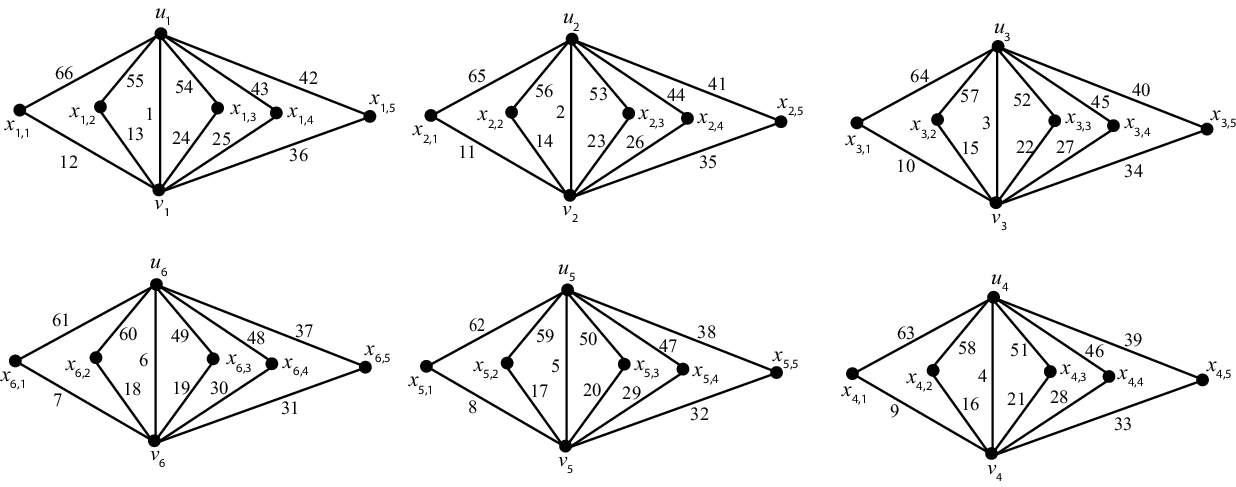, width=13cm}}
\caption{$6(P_2\vee O_5)$.}\label{fig:6P2VO5}
\end{figure}
\rsq
\end{example}

\begin{theorem}\label{thm-r(2sP2VO2n+1)} For $n\ge 1, r\ge 2, s\ge 1$, $\chi_{la}(r((2s)P_2 \vee O_{2n+1})) = 3$.
\end{theorem}

\begin{proof} Let $k=rs\ge 2$. Begin with $2k(P_2\vee O_{2n+1})$ and the local antimagic 3-coloring $f$ as defined in the proof of Theorem~\ref{thm-evenP2VO2n+1}. By Observation (4), we can now partition the $2k$ components into $2r$ blocks of $s$ component(s) with the $b$-th block $(1\le b\le 2r)$ containing vertices $u_i, v_i, x_{i,j}$ for $(b-1)s+1\le i\le bs, r\ge 2, s\ge 1$. Now, for each $j\in [1,2n+1]$, merge the vertices in $\{x_{i,j}, x_{2k+1-i, j} \;|\; (b-1)s + 1\le i\le bs\}$ to form a new vertex  of degree $4s$, still denoted $x_{i,j}$, $1\le i\le r$. We now get an $r$ components graph $r((2s)P_2 \vee O_{2n+1})$ with a bijective edge labeling such that the degree $2n+2$ vertices, $u_i$ and $v_i$, have induced vertex labels $(a) = f^+(u_i) = (n+1)(12kn+4k+1)+2k$ and $(b) = f^+(v_i) = (n+1)(4kn+4k+1)$ respectively. Moreover, the vertices $x_{i,j}$ have induced vertex label $(c) = f^+(x_{i,j}) = s(16kn+16k+2)$. Since $(a) > (b)$, we shall show that $(a), (b) \ne (c)$.

\ms\nt Now, $(a) - (c)$ is:
\begin{align}
&\quad \ 12kn^2 - 16kns + 16kn - 16ks + 6k + n - 2s+1\nonumber\\
& = 4kn(3n-4s+4) - 16ks + 6k + n - 2s + 1\label{eq:(a)-(c)-2}
\end{align}

\nt If $3n-4s+4\le 2$ (i.e., $-4s \le -3n-2$), then
\begin{align*}\eqref{eq:(a)-(c)-2} & \le 8kn - 16ks + 6k +n - 2s + 1\\
& \le  8kn-12kn-8k+6k+n-2s+1=-4kn-2k+n-2s+1<0.\end{align*}

\nt If $3n-4s+4\ge 3$ (i.e., $-4s\ge -3n - 1$), then
\begin{align*}\eqref{eq:(a)-(c)-2} & \ge 12kn - 16ks + 6k +n - 2s + 1\\
& \ge 2k+n-2s+1 =2rs+n-2s+1 > 0.\end{align*}

\nt Next, $(b)-(c)$ is :
\begin{align}
&\quad \ 4kn^2 - 16kns + 8kn - 16ks + 4k + n - 2s+1\nonumber\\
& = 4kn(n-4s+2)- 16ks + 4k + n-2s+1 \label{eq:(b)-(c)-2}
\end{align}
\nt If $n-4s+2\le 0$, then $\eqref{eq:(b)-(c)-2} \le -16ks+4k+2s-1 = -4k(3s-1)-s(4k-2)-1 < 0$.   If $n-4s+2 \ge 1$, then $\eqref{eq:(b)-(c)-2} \ge 4kn - 16ks + 4k + 2s \ge 4k(4s-1) - 16ks + 4k + 2s = 2s > 0$.

\ms\nt Therefore, $f$ is a local antimagic 3-coloring and $\chi_{la}(r((2s)P_2 \vee O_{2n+1})) = 3$. This completes the proof. \end{proof}

\nt Begin with $r((2s(P_2 \vee O_{2n+1})$ and the local antimagic 3-coloring $f$ defined in the proof of Theorem~\ref{thm-r(2sP2VO2n+1)}. In general, for $k = rs, n\ge 1, r\ge 2, s\ge 1$, we can also define $G_{r,2s}(2n+1)$ (similar to $G_{r,2s}(2n)$ as in Theorem~\ref{thm-G_{r,2s}(2n)}) obtained from $r((2s)P_2 \vee O_{2n+1})$ by splitting each $x_{i,j}$ of degree $4s$ into two vertices, say $y_{i,j}$ and $z_{i,j}$, of degree $2s$ such that both have equal induced vertex labels under $f$. By Observation (6), we know that $G_{r,2s}(2n+1)$ admits a bijective edge labeling, say $g$, such that $g^+(u_i) = f^+(u_i) = (n+1)(12kn+4k+1)+2k$, $g^+(v_i) = f^+(v_i) = (n+1)(4kn+4k+1)$ and $g^+(y_{i,j}) = g^+(z_{i,j}) = \frac{1}{2}f^+(x_{i,j}) = s(8kn+8k+1)$.

\begin{theorem}\label{thm-G_{r,2s}(2n+1)} For $n\ge 1$, $r\ge 2$, $s\ge 1$, $\chi_{la}(G_{r,2s}(2n+1)) = 3$. \end{theorem}

\begin{proof} By definition, we know that $k=rs\ge 2$ and $G_{r,2s}(2n+1)$ is an $r$ components bipartite graph with each component of equal partite sets size $2n+2s+1$. By  Lemma~\ref{lem-2part}, $\chi_{la}(G_{r,2s}(2n+1)) \ge 3$. Moreover, $G_{r,2s}(2n+1)$ admits a bijective edge labeling $g$ with induced vertex labels $(a) = g^+(u_i) = (n+1)(12kn+4k+1)+2k$, $(b) = g^+(v_i) =  (n+1)(4kn+4k+1)$, and $g^+(y_{i,j}) = g^+(z_{i,j}) =  s(8kn+8k+1)$. Since $(a) > (b)$, we shall show that $(a), (b) \ne (c)$.

\ms\nt Now, $(a) - (c)$ is:

\begin{align}
&\quad \ 12kn^2 - 8kns + 16kn - 8ks + 6k + n - s + 1\nonumber\\
& =  4kn(3n-2s+4) - 8ks + 6k + n - s +1      \label{eqn:(a)-(c)-2}
\end{align}
\nt If $3n-2s+4\le 2$ (i.e., $-2s \le -3n-2$), then
\begin{align*}
\eqref{eqn:(a)-(c)-2} & \le 8kn - 8ks + 6k + n -s +1 \\
 & \le 8kn - 12kn - 8k + 6k + n - s + 1 = -4kn - 2k + n - s + 1 < 0.
\end{align*}
\nt If $3n-2s+4 \ge 3$ (i.e., $-2s \ge -3n-1$), then
\begin{align*}
\eqref{eqn:(a)-(c)-2} & \ge 12kn - 8ks + 6k + n - s + 1 \\
 & \ge 2k+n-s+1 = 2rs + n - s + 1 > 0
\end{align*}

\ms\nt Next, $(b) - (c)$ is :
\begin{align}
& \quad \ 4kn^2 - 8kns + 8kn - 8ks + 4k + n - s + 1 \nonumber\\
& = 4kn(n - 2s+2)-8ks+4k+n-s+1 \label{eqn:(b)-(c)-2}
\end{align}
\nt If $n-2s+2\le 0$, then $\eqref{eqn:(b)-(c)-2} \le -8ks + 4k+s-1 = -2k(3s-2) - s(2k-1)-1 < 0$. If $n-2s+2\ge 1$, then $\eqref{eqn:(b)-(c)-2}\ge 4kn-8ks+4k+ s \ge 4k(2s-1) - 8ks + 4k + s = s > 0$.

\ms\nt Therefore, $g$ is a local antimagic 3-coloring and $\chi_{la}(G_{r,2s}(2n+1)) \le 3$. This completes the proof.               \end{proof}

\begin{example} Using the $6(P_2\vee O_5)$ as in Figure~\ref{fig:6P2VO5}, we get a $3(2P_2 \vee O_5)$ that admits a local antimagic 3-coloring with induced vertex labels $261, 111, 146$ respectively. Moreover, we can also get a $G_{3,2}(5)$ (see Figure~\ref{fig:G32(5)}) that admits a local antimagic 3-coloring with induced vertex labels $261, 111, 73$ respectively. 
\begin{figure}[H]
\centerline{\epsfig{file=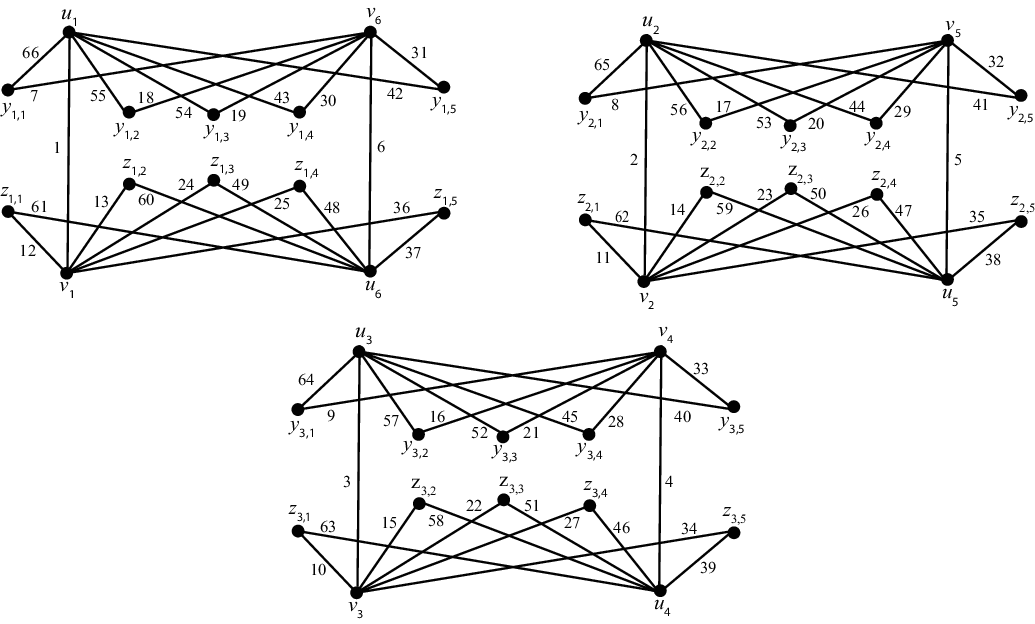, width=13cm}}
\caption{$G_{3,2}(5)$}\label{fig:G32(5)}
\end{figure}
 \rsq
\end{example}

\nt We can now define $\mathcal G_{r,2s}(2n+1)$ (similar to $\mathcal G_{r,2s}(2n)$ as in Corollary~\ref{cor-G_r,2s(2n)}) being the family of all non-isomorphic graphs that can be obtained by applying the delete-add process to $r((2s)P_2\vee O_{2n+1})$ and $G_{r,2s}(2n+1)$. The following corollary is also immediate.

\begin{corollary}\label{cor-G_r,2s(2n+1)} For $n\ge 1, r\ge 2, s\ge 1$, every graph $G\in \mathcal G_{r,2s}(2n+1)$ has $\chi_{la}(G) = 3$.  \end{corollary}

\nt Similar to Theorem~\ref{thm-G2n(2r+1,2s)}, we can also define $G_{2n+1}(2r+1,2s)$ and the following theorem. The proof is thus omitted.

\begin{theorem}\label{thm-G2n+1(2r+1,2s)} Suppose $n\ge 1$. If $2k=(2r+1)(2s)$ for $r,s\ge 1$, then $\chi_{la}(G_{2n+1}(2r+1,2s)) = 3$. \end{theorem}

 \nt Similar to Corollary~\ref{cor-G2n(2r+1,2s)}, we also have the following corollary immediately.

\begin{corollary} For $n,r,s\ge 1$, every graph $G \in \mathcal G_{2n+1}(2r+1,2s)$ has $\chi_{la}(G) = 3$. \end{corollary}

\nt Similar to Theorems~\ref{thm-J1} and~\ref{thm-J2}, we also let $\mathcal J_1(k,2n+1,s)$ (respectively, $\mathcal J_2(k,2n+1,s)$) be the family of graphs obtained from $k(2P_2\vee O_{2n+1})$ (respectively, $G_{k,2}(2n+1)$) after partitioning the $2k$ vertices $u_1, \ldots, u_{2k}$, with induced vertex label $ (n+1)(12kn+4k+1)+2k$.

\begin{theorem}\label{thm-J1-2} For $k,s\ge 2$, $n\ge 1$, each graph $G\in \mathcal J_1(k,2n+1,s)$ has $\chi_{la}(G) = 3$. \end{theorem}

\begin{proof} By definition, $\chi_{la}(G) \ge \chi(G) = 3$. Keeping the bijective edge labeling defined for $k(2P_2\vee O_{2n+1})$, we can immediately conclude that $G$ admits a bijective edge labeling with induced vertex labels $(a) = s[(n+1)(12kn+4k+1) + 2k]$, $(b) = (n+1)(4kn+4k+1)$ and $(c) = 16kn+16k+2$. It is easy to show that $(a) > (b) > (c)$ and the labeling is local antimagic. Thus, $\chi_{la}(G) \le 3$. Consequently, $\chi_{la}(G) = 3$.
\end{proof}

\nt Similarly, we also have the following theorem with the proof omitted.

\begin{theorem}\label{thm-J2-2} Suppose $k,s\ge 2, n\ge 1$ and $G\in \mathcal J_2(k,2n,s)$. If $G$ is bipartite with every component of equal partite set size or $G$ is tripartite, then $\chi_{la}(G) = 3$. \end{theorem}

\nt Similar to Corollaries~\ref{cor-H1} and~\ref{cor-H2}, we also define $t\ge 2$ components graphs $H_1(k_1,2n+1) + H_1(k_2,2n+1) + \cdots + H_1(k_t,2n+1)$ and $H_2(k_1,2n+1) + H_2(k_2,2n+1) + \cdots + H_2(k_t,2n+1)$ using $k(P_2\vee O_{2n+1})$ and $G_{k,2}(2n+1)$ respectively.

\begin{corollary} For $n\ge 1$ and $t, k_a\ge 2$, the graph $G_1 = H_1(k_1,2n+1) + H_1(k_2,2n+1) + \cdots + H_1(k_t,2n+1)$ and $G_2 = H_2(k_1,2n+1) + H_2(k_2,2n+1) + \cdots + H_2(k_t,2n+1)$ has $\chi_{la}(G_i) = 3$ for $i=1,2$. \end{corollary}

\nt Further note that the following graphs are $(2n+2)$-regular:
\begin{enumerate}[(1)]
\item $(n+1)P_2\vee O_{2n+1}$.
\item $r((n+1)P_2\vee O_{2n+1})$.
\item $G_{r,2n+2}(2n+1)$.
\item each graph in $\mathcal G_{r,2n+2}(2n+1)$.
\end{enumerate}

\nt Thus, we have the following corollary.

\begin{corollary} For each $n\ge 1$, there are (possibly connected) $(2n+2)$-regular bipartite or tripartite graphs with local antimagic chromatic number $3$. \end{corollary}

\section{Conclusion and Open Problems}

\nt This paper is a natural extension to~\cite{LauShiu-join, LauShiu-chila3}. We successfully proved that the join of a  1-regular graph and a null graph has local antimagic chromatic number 3. Moreover, $\chi_{la}(a(bP_2\vee O_m)) = 3$ for all $a, b, m\ge 2$ except $a$ is even and $b$ is odd.

\begin{conjecture} For $r, s\ge 1$ and $m\ge 2$, $\chi_{la}(2r((2s+1)P_2 \vee O_m)) = 3$.  \end{conjecture}

\nt Interested readers can refer to~\cite{LSN-DMGT, LPAS-DMAA} for more results on local antimagic chromatic number of the join of a path or a cycle and a null graph. We end this paper with the following problem.

\begin{problem} Let $G$ be a path or a cycle. Determine $\chi_{la}(mG \vee O_n)$ for $m, n\ge 2$.  \end{problem}

\end{document}